\documentclass[10pt,reqno,oneside]{amsart}
\usepackage{amsmath,amsthm,amsfonts,amssymb}
\usepackage{indentfirst}
\usepackage{url}
\usepackage{hyperref}
\usepackage{graphicx}
\usepackage{color}

\theoremstyle{plain}
\newtheorem{theorem}{Theorem}[section]

\newtheorem{prop}[theorem]{Proposition}
\theoremstyle{definition}
\newtheorem{defn}[theorem]{Definition}
\newtheorem{exa}[theorem]{Example}
\newtheorem{obs}[theorem]{Remark}
\newcommand{\fim}{{$\hfill\square$\vspace{.2cm}}}

\begin{document}
	
	\baselineskip=17pt

	\title[The GMS model with threshold extinction]{The GMS model with threshold extinction}
	\author[Luiz Renato Gon\c calves Fontes]{Luiz Renato Fontes}
	\author[Carolina Grejo]{Carolina Grejo}
		\author[F\'abio Sternieri Marques]{F\'abio Sternieri Marques}
		\address[L. R. Fontes and  F. S Marques]{USP - Instituto de Matem\'atica e Estat\'istica, 
			Rua do Mat\~ao, 1010 Cidade Universit\'aria
			Butant\~a, 05508090 - S\~ao Paulo, SP - Brasil, SP, Brasil}
		\email{lrfontes@usp.br, fmarques2015@gmail.com}
	\address[C. Grejo ]{UFSCar - Departamento de Estat\'istica, 
		Rodovia Washington Luiz, km 235, CEP 13565-905, São Carlos, SP, Brazil}

	\email{carolina@ime.usp.br}
	
\thanks{Research supported by CAPES; PNPD-CAPES; CNPq 311257/2014-3; FAPESP 2017/10555-0}

	\keywords{Poisson Processes, Records, Evolution.}
	\subjclass[2010]{60J20, 60K35}
	\date{\today}

	\begin{abstract} 
	We propose a variation of the GMS model of evolution of species. In this version, as in the GMS model, at each birth, the new species in the system is labeled with a random fitness mark, but in our variation, to each extinction event is associated a random threshold mark and all species with fitness lower than the threshold are removed from the system. 
	We present necessary and suficient criteria for the recurrence and transience of the empty configuration of species; we show the existence of a long time limit distribution of species in the system, and present necessary and suficient criteria for the finiteness of the number of species in that distribution. There is a remarkable symmetry between both sets of criteria.
	\end{abstract}
	
\maketitle
	\section{Introduction}
	\label{S: Introduction}
The GMS model, first proposed by Guiol et al in \cite{GMS11}, describes the evolution (in discrete time) of species who independently appear at each time step with a given fixed probability, and are assigned a {\em fitness} random variable, with a fixed continuous distribution. Extinction occurs at each time step, also independently and with a fixed probability, whenever there is at least one species present at the corresponding time, in which case the one with the least fitness gets extinct. Several variations of this model were studied, as in Ben Ari et al \cite{IDO11},  Skevi and Volkov \cite{SV12}, Bertacchi et al \cite{BLZ18} and Grejo at al \cite{GMR16}.  Further, Guiol et al \cite{GMS13} proposed a variation for the model, where the evolution is given in continuous time. See also \cite{FS16} for a closely related model with a different motivation. We refer to the above literature for further motivation of the models treated therein and results. 

We consider here a variation of the GMS model \cite{GMS13} in which, as in that model, new species are born at a given rate, and at possibly another rate we have extinction of species. For each new species in the system we associate a positive random number, chosen from a distribution $F_{\star}$. We call this random number the {\em fitness} of the species. So far, the setting is the same as in  \cite{GMS13}. 

Our variation is related with the extinction events. At the time of each such event, we have a positive threshold random variable, with distribution $F_{\dagger}$, and all species with fitness below this threshold at that time get extinct. 

We believe this is a natural variation of the original model, when we consider extinction of species in the natural world produced by major events such as abrupt habitat change, where conceivably each species is affected according to its own aptitude to face the new challenge,  irrespective of other species.

The random fitnesses and the random thresholds are  independent of each other and of everything else in the process. We assume $F_{\star}$ and $F_{\dagger}$ are  continuous on $\mathbb{R}_+=[0,\infty)$. So, in particular, we can and will identify each species with its fitness. We also assume, for simplicity, that $F_{\star}$ and $F_{\dagger}$ have unbounded support.

Let ${\Pi}_{\star}$ and ${\Pi}_{\dagger}$ be independent Poisson point processes with rates $\lambda_{\star}$ and $\lambda_{\dagger}$, respectively. Define $\{T_i\}_{i\in {\mathbb{Z}}^*}$ as the set of birth time instants of a new species in the system, define $\{S_j\}_{j\in {\mathbb{Z}}^*}$ as the instants of time in which there is an extinction event. These sequences represent the points in the Poisson processes and are indexed in increasing order. At each time $T_i$ we assign to 
the newly appeared species the fitness $X_i$, drawn from $F_{\star}$, and to each time  $S_j$ we associate the threshold $Y_j$ from distribution $F_{\dagger}$.

Given a locally finite subset $A$ of ${\mathbb{R}}_+$,
let $\eta_t=\{\eta_t(s),~s \geq t\}$ be the GMS model with threshold extinction starting from $A$ at time $t$.

So, at time $t\in \mathbb{R}$, the process has initial configuration $\eta_t(t)=A,$ 
and at time $s\in (t,\infty)$, the process has the configuration $\eta_t(s)$ which is composed of all species/fitnesses either of $A$ or that have appeared in the time interval $(t,s]$, and that have survived the events of extinction in $[t,s]$, i. e.  species whose fitnesses are greater than the highest threshold drawn in events of extinction in $[t,s]$
after their birth. 

Notice that as long as $A$ does not depend on $t$, neither does the distribution of $\eta_t$, so below we will often restrict to $\eta_0$.

\section{Results}

We derive three kinds of results. First, criteria for recurrence and transience of the empty configuration in $\eta_0$; see Theorem~\ref{teo1} below; they are obtained from the analysis of a Poisson process of records. Second,
we derive the existence of a limit distribution for $\eta_0(s)$ as $s\to\infty$; see Theorem~\ref{teo2} below. Finally, we derive criteria for finitude and infinitude of the number of species present in the limit distribution; see Theorem~\ref{teo3} below; curiously, the Poisson process of records of the recurrence and transience issue appears here as well, but {\em in reverse}. 
The behavior of $\bar{F_{\dagger}}\circ \bar{F_{\star}}^{-1}$ at the origin plays a determinant role in the first result,
and thus, in reverse, so does that of $\bar{F_{\star}}\circ \bar{F_{\dagger}}^{-1}$ in the third one; 
as usual, for $\ast=\star$ and $\dagger$,
$\bar{F_{\ast}}=1-F_\ast$, and 
$\bar{F_{\ast}}^{-1}$ indicates the (right-continuous) inverse of $\bar{F_{\ast}}$. 
Proofs are deferred to the last section.

Let us define the functions  $R_{\star}, R_{\dagger}:\mathbb{R}_+\rightarrow\mathbb{R}_+$ by
\begin{align*}
R_{\star}(x)&=-\log{\bar{F_{\star}}(x)},\\
R_{\dagger}(x)&=-\log{\bar{F_{\dagger}}(x)}.
\end{align*}

For birth events denote by $I_k$ the record indexes  and by $X_{I_k}$ the record values,  as follows: $I_1\doteq1$; and 
for $k\geq1$
\begin{equation*}
I_{k+1}\doteq\min\{i>I_k:X_i>X_{I_k}\}. 
\end{equation*}

\begin{prop} [Proposition 4.11.1, Resnick \cite{R73}]
	\label{rec}
	The record values $(X_{I_k})_{k\geq1}$ form  a  Poisson point process in $\mathbb{R}_+$ with intensity measure $\int_B R_{\star}(dx)$, $B \in \mathcal{B}(\mathbb{R}_+)$.
\end{prop}

\begin{prop}\label{marc:lambda}
	Let $T_{I_k}$ be the time of the  $k-th$ record and denote by $\Delta T_{I_k}\doteq T_{I_{k+1}}-T_{I_k}$ the interval between two consecutive records. Then, $\{(X_{I_k},\Delta T_{I_k})_{k\geq1}\}$ is a Poisson process in $\mathbb{R}_+^2$ with intensity measure
	\begin{equation*}
	\hat{\mu}_{\star}(C)\doteq\iint_C \, \lambda_{\star}\bar{F_{\star}}(x)e^{-\lambda_{\star}\bar{F_{\star}}(x)s}ds R_{\star}(dx), \text{}\quad C\in\mathcal{B}(\mathbb{R}_+^2).
	\end{equation*}
	\end{prop}
	
\begin{prop}\label{marc:mu}
	The set of points $(S_j,Y_j)_{j\geq1}$ is a Poisson process in  $\mathbb{R}_+^2$, denoted by $\hat{\Pi}_{\dagger}$, with intensity measure
	\begin{equation*}
	\mu_{\dagger}(C)\doteq\lambda_{\dagger}\iint_C dt\, F_{\dagger}(dx), \text{}\quad C\in\mathcal{B}(\mathbb{R}_+^2).
	\end{equation*}
\end{prop}
To study recurrence and transience, we will build a ladder of records using the process of births and the fitness associated to each species.
Let us denote the random region above each step of the ladder by $\{D_k\}_{k\geq 1}$, and the region above the full ladder is denoted by $D$, namely, for $k\geq1$
\begin{equation*}
D_k\doteq [T_{I_k},T_{I_{k+1}})\times [X_{I_k},\infty) \quad \text{and} \quad D\doteq \bigcup_{k\geq1} D_k.
\end{equation*}

 \begin{equation*}
	\Lambda\doteq\mu_{\dagger}(D);\quad
	M\doteq\#(D\cap\hat{\Pi}_{\dagger}).
\end{equation*}

We may use Campbell's formula to find
\begin{equation}\label{eq:expM}
	\mathbb{E}[M]= 
	\frac{\lambda_{\dagger}}{\lambda_{\star}}\int_0^\infty\frac{\bar{F_{\dagger}}(x)}{\bar{F_{\star}}(x)}R_{\star}(dx).
\end{equation}

 Observe that  $M$ denotes the number of extinction events in $D$ and, given $D$, it has a Poisson distribution with mean $\Lambda$ (because $\hat{\Pi}_{\dagger}$ is a Poisson process), so 

 \begin{equation}\label{eq.exp}
 \mathbb{E}[e^{-tM}]=\mathbb{E}\big[\mathbb{E}[e^{-tM}|\Lambda]\big]=\mathbb{E}\big[e^{-(1-e^{-t})\Lambda}\big].
 \end{equation}

Here we allow $\Lambda=\infty$, in which case $M=\infty$ a.s.

 \begin{obs}
 	Define $h_{\dagger}: \mathbb{R}_+^2 \rightarrow \mathbb{R}_+$ by
 $h_{\dagger}(x,s)\doteq \lambda_{\dagger}s\bar{F_{\dagger}}(x)$.
 	From the definition of $\mu_{\dagger}$ in Proposition \ref{marc:mu}, we have
 	
 	\begin{equation}\label{eq.Delta}
 	\Lambda=\mu_{\dagger}(D)=\sum_{k\geq1}\mu_{\dagger}(D_k)=\sum_{k\geq1}\lambda_{\dagger}\Delta T_{I_k} \bar{F_{\dagger}}(X_{I_k})=\sum_{k\geq1}h_{\dagger}(X_{I_k},\Delta T_{I_k}).
 	\end{equation}
 \end{obs}

\begin{prop}\label{LaplaceM2} 
	For $t>0$
	\begin{equation*}
	\label{eq:espM}
	\mathbb{E}[e^{-tM}]=\exp{\Bigg[ -\int_0^\infty \frac{(1-e^{-t}) \lambda_{\dagger} \bar{F_{\dagger}}(x)}{\lambda_{\star}\bar{F_{\star}}(x)+(1-e^{-t}) \lambda_{\dagger} \bar{F_{\dagger}}(x)} R_{\star}(dx)\Bigg]}.
	\end{equation*}
\end{prop}

\subsection{Recurrence/transience of $\eta_0$}\label{transrec}

We start by discussing the a.s.~finitude of $M$.
Define the functions $\phi:\mathbb{R}_+\rightarrow\mathbb{R}_+$ and $\psi:[0,1]\rightarrow[0,1]$ 
by
\begin{equation*}\label{eq.lambda}
\phi(t)\doteq \int_0^\infty\frac{(1-e^{-t}) \lambda_{\dagger} \bar{F_{\dagger}}(x)}{\lambda_{\star}\bar{F_{\star}}(x)+(1-e^{-t}) \lambda_{\dagger} \bar{F_{\dagger}}(x)}R_{\star}(dx),
\end{equation*}
and
\begin{equation}\label{eq:comp}
\psi(u)\doteq \bar{F_{\dagger}}\circ \bar{F_{\star}}^{-1}(u).
\end{equation}

By monotone convergence, we have that
\begin{equation*}
\phi(\infty)\doteq \lim_{t\rightarrow \infty}\phi(t)=\int_0^\infty\frac{\lambda_{\dagger} \bar{F_{\dagger}}(x)}{\lambda_{\star}\bar{F_{\star}}(x)+\lambda_{\dagger} \bar{F_{\dagger}}(x)}R_{\star}(dx).
\end{equation*}

\begin{prop}
	\label{prop:equiv}
		The following statements are equivalent:
		\renewcommand{\labelenumi}{\roman{enumi})}
		\begin{enumerate}
			\item $\phi(\infty)<\infty$;
			\item $\mathbb{P}(M<\infty)=1$;
			\item $\mathbb{E}[M]<\infty$;
			\item $\dfrac{\psi(u)}{u^2}$ is integrable at the origin;
			\item $\psi(s^{-1})$ is integrable at infinity.
		\end{enumerate}
\end{prop}

\begin{defn}
	The process $\eta_0$ is said to be transient if $\{s\geq0:\eta_0(s)=\varnothing\}$ is a.s.~a bounded set.
	On the other hand, if $\{s\geq0: \eta_0(s)=\varnothing\}$ 
	 is a.s.~an unbounded set, we say that $\eta_0$ is recurrent.
\end{defn}

	\begin{theorem}\label{teo1}
	The process $\eta_0$ is transient if 
	\begin{equation*}
	\int_0^\infty\frac{\bar{F_{\dagger}}(x)}{\bar{F_{\star}}(x)} R_{\star}(dx)<\infty,
	\end{equation*}
	and recurrent otherwise. 
\end{theorem}

\begin{exa}
	\label{ex:rec/trans}
	Suppose that $(X_i)$ and $(Y_j)$ are exponentially distributed with parameter ${\alpha}_{\star}$ and ${\alpha}_{\dagger}$, respectively. Then,
	\begin{equation*}
	\int_0^\infty\frac{\bar{F_{\dagger}}(x)}{\bar{F_{\star}}(x)} R_{\star}(dx)=\left\{ 
	\begin{array}{cc}
	\frac{{\alpha}_{\star}}{{\alpha}_{\dagger}-{\alpha}_{\star}}, & \mbox{if}~{\alpha}_{\dagger}>{\alpha}_{\star}\\
	\infty, & \mbox{if}~{\alpha}_{\dagger}\leq{\alpha}_{\star}.
	 \end{array} 
	 \right.
	\end{equation*}

By Proposition \ref{prop:equiv}, if ${\alpha}_{\dagger}\leq{\alpha}_{\star}$, then $M=\infty$ a.s.,~and from Theorem \ref{teo1}, $\eta_0$ is recurrent. If ${\alpha}_{\dagger}>{\alpha}_{\star}$, then $\eta_0$ is transient. We may compute the distribution of $M$ in this case (it will come in handy below), as follows:
$\mathbb{E}\left[ e^{-t\Lambda}\right]=\left( 1+t\beta^{-1}\right)^{-r}$,
and, by Proposition \ref{eq:espM}, 
$\mathbb{E}[e^{-tM}]=\left( \frac{1-p}{1-pe^{-t}}\right)^r,$
where $r=\frac{{\alpha}_{\star}}{{\alpha}_{\dagger}-{\alpha}_{\star}},~\beta=\frac{{\lambda}_{\star}}{{\lambda}_{\dagger}}$ and $p=\frac{{\lambda}_{\dagger}}{{\lambda}_{\star}+{\lambda}_{\dagger}}$.
Hence, $\Lambda$ follows the gamma distribution with parameters $r$ and $\beta$, and $M$ follows the negative binomial distribution with parameters $r$ and $p$.
\end{exa}

\subsection{Existence and in/finitude of a long time limit distribution} We now establish the existence of a limit distribution for ${\eta}_0(s)$ as $s\to \infty$. Remarkably, a ladder construction based on a record process comes up here as well, entirely parallel to that of Subsection~\ref{transrec}, with births and extinctions swapping roles. This immediately yields necessary and sufficient criteria for the almost sure in/finitude of the number of species present in the limit distribution, identical to those for transience/recurrence of $\eta_0$, except that the symbols $\star$ and $\dagger$ swap roles.

As a preliminary, we will use the process of extinctions and the associated thresholds to build a ladder of records, much as in Subsection~\ref{transrec}.
We define the k-$th$ record index $J_k$ and the record value $Y_{J_k }$ 
as follows: $J_1\doteq-1$ and; for $k\geq1$
\begin{equation*}\label{ind}
J_{k+1}\doteq\max\{j<J_k:Y_j>Y_{J_k}\}
\end{equation*}

Again by Proposition 4.11.1 in~\cite{R73}, we get that
the 
$\{Y_{J_k}\}_{k\geq1}$ form a Poisson point process in $\mathbb{R}_+$ with intensity measure 
$\int_B R_{\dagger}(dx)$,  $B \in \mathcal{B}(\mathbb{R}_+).$

Denote by  $S_{J_k}(0)$ the time of the $k$-th record, and by 
$\Delta S_{J_k}\doteq S_{J_k}-S_{J_{k+1}}$ the time span between two consecutive records. Similarly as in Subsection~\ref{transrec}, we have the following results.

\begin{prop}\label{marc}
The points $(Y_{J_k},\Delta S_{J_k})_{k\geq1}$ form a Poisson point process in $\mathbb{R}_+^2$ with intensity
	\begin{equation*}
	\hat{\mu}_{\dagger}(C)\doteq\iint_C\, \lambda_{\dagger}\bar{F_{\dagger}}(x)e^{-\lambda_{\dagger}\bar{F_{\dagger}}(x)s}ds  R_{\dagger}(dx),\,\,\, C\in\mathcal{B}(\mathbb{R}_+^2).
	\end{equation*}
\end{prop}

\begin{prop}\label{marc:lamb}
	The points $\{(T_{-i},X_{-i})\}_{i\geq1}$ form a Poisson point process in $\mathbb{R}_-\times\mathbb{R}_+$,
	$\mathbb{R}_-=(-\infty,0]$,  denoted by $\hat{\Pi}_{\star}$, with intensity
	\begin{equation*}
	\mu_{\star}(C')\doteq\lambda_{\star}\iint_{C'}dt\, F_{\star}(dy),\,\,\, C'\in\mathcal{B}(\mathbb{R}_-\times\mathbb{R}_+).
	\end{equation*}
\end{prop}
We thus have our ladder of thresholds; denote the region above each step by
\begin{equation*}
E_0\doteq[S_{J_1},0)\times\mathbb{R}_+;
\end{equation*}
\begin{equation*}
E_k\doteq[S_{J_{k+1}},S_{J_k})\times[Y_{J_k},\infty),\,k\geq1; \,\,\,\,
E\doteq\bigcup_{k\geq1}E_k.
\end{equation*}

We state our existence result.
\begin{theorem}\label{teo2}
	$\eta_0(t)$ converges in distribution to  $\hat{\eta}$ as $t\rightarrow\infty$, where 
	\begin{equation*}
	\hat{\eta}\doteq\{ X_i:T_i \in (S_{J_1},0]\}\cup\Big(\bigcup_{k\geq1}\left\{X_i>Y_{J_k}:T_i \in \left(S_{J_{k+1}},S_{J_k}\right]\right\}\Big) .
	\end{equation*}
\end{theorem}

\begin{obs}
	The topology for weak convergence is the usual one in the context of point processes. Our proof indeed makes use of a coupling to a sequence of processes for which the convergence is a strong one, and follows by monotonicity. 
\end{obs}

Next, we address the issue of finitude of the number of species in $\hat{\eta}$.
For each $k\geq0$, let
\begin{equation*}
\Sigma_k\doteq\mu_{\star}(E_k); \quad
N_k\doteq\#\{E_k\cap \hat{\Pi}_{\star}\}. 
\end{equation*}

Let also
$\Sigma\doteq\mu_{\star}(E)$;
$N\doteq\#\{E\cap \hat{\Pi}_{\star}\}.$
We may use Campbell's formula to find
\begin{equation*}
	\mathbb{E}[N]=\frac{\lambda_{\star}}{\lambda_{\dagger}}\int_0^\infty\frac{\bar{F_{\star}}(x)}{\bar{F_{\dagger}}(x)}R_{\dagger}(dx).
\end{equation*}

Note that $N$ is the number of birth events above the threshold ladder. Also, note that  $E$  has a parallel structure to that of the $ D $ ladder of Subsection~\ref{transrec}; the random variables $\Sigma$ and $N$ are parallel to $\Lambda$ and $M$ in that same subsection. Thus, we get parallel results, once we exchange the roles of $({\lambda}_{\star}, {\bar{F}}_{\star})$ and $({\lambda}_{\dagger},{\bar{F}}_{\dagger})$.

From Theorem \ref{teo2}, the number of species present in the limit distribution $\hat{\eta}$, denoted by $\#\hat{\eta}$,  is
\begin{equation}\label{numlim}
\#\hat{\eta}=\sum_{k\geq0}\#\{E_k\cap \hat{\Pi}_{\star}\}=\sum_{k\geq0}N_k=N_0+N.
\end{equation}

It is enough to consider the finitude of $N$. From the parallel situation of Subsection~\ref{transrec}, we get the following results. For $t>0$
 \begin{eqnarray*}\nonumber
	&\mathbb{E}[e^{-tN}]=\mathbb{E}[e^{-(1-e^{-t})\Sigma}]&\\
&=\exp{\Bigg[-\displaystyle\int_0^\infty\frac{(1-e^{-t})\lambda_{\star}\bar{F_{\star}}(x)}{\lambda_{\dagger}\bar{F_{\dagger}}(x)+(1-e^{-t})\lambda_{\star}\bar{F_{\star}}(x)}R_{\dagger}(dx)\Bigg]}.&
\end{eqnarray*}

Setting 
\begin{equation*}
\bar{\phi}(t)\doteq \int_0^\infty\frac{(1-e^{-t})\lambda_{\star}\bar{F_{\star}}(x)}{\lambda_{\dagger}\bar{F_{\dagger}}(x)+(1-e^{-t})\lambda_{\star}\bar{F_{\star}}(x)}R_{\dagger}(dx),\quad t>0,
\end{equation*}
and
\begin{equation*}
\bar{\psi}(u)\doteq \bar{F_{\star}}\circ \bar{F_{\dagger}}^{-1}(u),\, u>0,
\end{equation*}

we have
\begin{equation*}
\bar{\phi}(\infty)\doteq\lim_{t\rightarrow\infty}\bar{\phi}(t)\\
=\int_0^\infty\frac{\lambda_{\star} \bar{F_{\star}}(x)}{\lambda_{\dagger}\bar{F_{\dagger}}(x)+\lambda_{\star} \bar{F_{\star}}(x)}R_{\dagger}(dx).
\end{equation*}

\begin{prop}\label{equiv2}
	The following statements are equivalent:
	\renewcommand{\labelenumi}{\roman{enumi})}
	\begin{enumerate}
		\item \label{i1} $\bar{\phi}(\infty)<\infty$;
		\item \label{i2} $\mathbb{P}(N<\infty)=1$;
		\item \label{i3} $\mathbb{E}[N]<\infty$.
		\item \label{i4} $\dfrac{\bar{\psi}(u)}{u^2}$ is integrable at the origin;
		\item \label{i5} $\bar{\psi}(s^{-1})$ is integrable at infinity.
	\end{enumerate}
\end{prop}
\begin{theorem}\label{teo3}
	The number of species in the limit distribution, $\#\hat{\eta}$, is finite if
	\begin{equation*}
	\int_0^\infty\frac{\bar{F_{\star}}(x)}{\bar{F_{\dagger}}(x)} R_{\dagger}(dx)<\infty,
	\end{equation*}
	and infinite otherwise. 
\end{theorem}

\begin{exa}
	Let $(X_i)$ and $(Y_j)$ be as in Example \ref{ex:rec/trans}. We have then
\begin{equation*}
	\int_0^\infty\frac{\bar{F_{\star}}(x)}{\bar{F_{\dagger}}(x)} R_{\dagger}(dx)=\left\{ 
	\begin{array}{cc}
	\frac{{\alpha}_{\dagger}}{{\alpha}_{\star}-{\alpha}_{\dagger}}, & \mbox{if}~{\alpha}_{\star}>{\alpha}_{\dagger}\\
	\infty, & \mbox{if}~{\alpha}_{\star}\leq{\alpha}_{\dagger}.
	\end{array} 
	\right.
\end{equation*}

By the Theorem~\ref{teo3}, if ${\alpha}_{\star}\leq{\alpha}_{\dagger}$, then $\#\hat{\eta}=\infty$ a.s. If ${\alpha}_{\star}>{\alpha}_{\dagger}$, then $\#\hat{\eta}<\infty$ a.s.; let us find its distribution. We can show that, as in Example~\ref{ex:rec/trans}, $N$ follows the negative binomial distribution with parameters $\frac{{\alpha}_{\dagger}}{{\alpha}_{\star}-{\alpha}_{\dagger}}$ and $\frac{{\lambda}_{\star}}{{\lambda}_{\star}+{\lambda}_{\dagger}}$. 
From the joint distribution of $(N_0,\Sigma_0)$, discussed in the proof of Theorem~\ref{teo3} below, we have that $N_0$  follows the negative binomial distribution with parameters $1$ and $\frac{{\lambda}_{\star}}{{\lambda}_{\star}+{\lambda}_{\dagger}}$. 
It may be checked that $N_0$ and $N$ are independent, and we find from
(\ref{numlim}) that $\#\hat{\eta}$ follows the negative binomial distribution with parameters 
$\frac{{\alpha}_{\star}}{{\alpha}_{\star}-{\alpha}_{\dagger}}$
and 
$\frac{{\lambda}_{\star}}{{\lambda}_{\star}+{\lambda}_{\dagger}}$.
\end{exa}

\begin{obs}
	As final remark, we might term the case of $\infty$ in both criteria in Theorems \ref{teo1} and \ref{teo3} as {\em null recurrent}. This is of course the case when $F_{\star}=F_{\dagger}$, which then becomes a natural candidate for comparison with the GMS model {\em below the critical point}, which also is recurrent and has infinitely many species in its limiting distribution.
\end{obs}

\section{Proofs}
Our proofs of Propositions \ref{marc:lambda} and \ref{marc:mu} use Definition 5.3 (K-marking process), as well as Theorem 5.6 (Marking Theorem) in \cite{LP17} (on pages 40 and 42, respectively).

\medskip

\noindent\textbf{\textit{Proof of Proposition \ref{marc:lambda}.}}~Conditional on $\{(X_{I_k})_{k\geq1}=(x_k)_{k\geq1}\}$, the random variables $\Delta T_{I_k}$ are independent of each other and follow the exponential distribution of rate $\lambda_{\star}\bar{F_{\star}}(x_k)$. Denote by $K_{\star}(x,B)$ the following probability kernel: for $x\geq0$, and $B=(s_1,s_2]\subset[0,\infty)$,  let
$K_{\star}(x,(s_1,s_2])\doteq e^{-\lambda_{\star}\bar{F_{\star}}(x)s_1}-e^{-\lambda_{\star}\bar{F_{\star}}(x)s_2}$, 
and extend the definition for Borelians $B$ in the usual way.
 
The points $(X_{I_k},\Delta T_{I_k})_{k\geq1}$ form a \mbox{$K_{\star}$-marking} of the Poisson process of Proposition \ref{rec}. The result follows by the Marking Theorem  (Theorem 5.6 of \cite{LP17}).
\fim

\noindent\textbf{\textit{Proof of Proposition \ref{marc:mu}.}}~
The random variables $(Y_j)_{j\geq1}$ are independent of each other and independent of 
$(S_j)_{j\geq1}$. Then, the points $(S_j,Y_j)_{j\geq1}$ form a  {$\mathbb{Q}_{\dagger}$-marking independent of $(S_j)_{j\geq1}$, where $\mathbb{Q}_{\dagger}$ is the measure induced by $Y$. The result follows again by the Marking Theorem.
\fim

\noindent\textbf{\textit{Proof of Proposition \ref{LaplaceM2}.}}
From equations (\ref{eq.exp}) and (\ref{eq.Delta})
\begin{eqnarray}
&&\mathbb{E}[e^{-tM}]=\mathbb{E}[e^{-(1-e^{-t})\Lambda}]=\mathbb{E}[e^{-\sum_{k\geq1}(1-e^{-t})h_{\dagger}(X_{I_k},\Delta T_{I_k})}]\nonumber\\
&&=\exp{\Bigg[-\iint_{[0,\infty)^2} (1-e^{-(1-e^{-t})h_{\dagger}(x,s)})\hat{\mu}_{\star}(dx,ds)\Bigg] } \label{FuncLap} \\
&&=\exp{\Bigg[-\iint_{[0,\infty)^2}(1-e^{-(1-e^{-t})\lambda_{\dagger}s\bar{F_{\dagger}}(x)})\lambda_{\star}\bar{F_{\star}}(x) e^{-\lambda_{\star}\bar{F_{\star}}(x)s}ds R_{\star}(dx)\Bigg]}. \label{eq.prop}
\end{eqnarray}

In (\ref{FuncLap}) we used the characterisation of a Poisson process via its Laplace functional (see Theorem 3.9 on page 23 of \cite{LP17}), and (\ref{eq.prop}) follows by Proposition \ref{marc:lambda}.

Integrating the above expression in $s$, we have
\begin{equation*}
\begin{split}
\mathbb{E}[e^{-tM}]&=\exp{\bigg[-\int_0^\infty\bigg(1-\frac{\lambda_{\star}\bar{F_{\star}}(x)}{\lambda_{\star}\bar{F_{\star}}(x)+(1-e^{-t}) \lambda_{\dagger} \bar{F_{\dagger}}(x)}\bigg) R_{\star}(dx)\bigg]}\\
&=\exp{\bigg[ -\int_0^\infty \frac{(1-e^{-t}) \lambda_{\dagger} \bar{F_{\dagger}}(x)}{\lambda_{\star}\bar{F_{\star}}(x)+(1-e^{-t}) \lambda_{\dagger} \bar{F_{\dagger}}(x)} R_{\star}(dx)\bigg]}. \hspace{3,2cm}\square
\end{split}
\end{equation*}
\vspace{.1cm}

\noindent\textbf{\textit{Proof of Proposition \ref{prop:equiv}.}} 
(i)$\Leftrightarrow$(ii) and (iv)$\Leftrightarrow$(v) are straightforward.
Let us show that (iii)$\Rightarrow$(i). From
\begin{equation*}
\frac{\lambda_{\dagger} \bar{F_{\dagger}}(x)}{\lambda_{\star}\bar{F_{\star}}(x)+\lambda_{\dagger} \bar{F_{\dagger}}(x)}\leq\frac{\lambda_{\dagger} \bar{F_{\dagger}}(x)}{\lambda_{\star}\bar{F_{\star}}(x)},
\end{equation*}
and integrating against $R_{\star}(dx)$, we get $\phi(\infty)\leq\mathbb{E}[M].$
Thus, (iii)$\Rightarrow$(i).

To show that (i)$\Rightarrow$(iii), we will change variables.
Let $y=R_{\star}(x)=-\log \bar{F_{\star}}(x)$, so 
$x=\bar{F_{\star}}^{-1}(e^{-y})$.
From this and (\ref{eq:comp}),
$\bar{F_{\dagger}}(x)=\bar{F_{\dagger}}\circ \bar{F_{\star}}^{-1}(e^{-y})=\psi(e^{-y}).$
Making $u=e^{-y}$, we find that 
\begin{equation*}\label{troca}
\begin{split}
\phi(\infty)&=\int_0^\infty\frac{\lambda_{\dagger} \bar{F_{\dagger}}(x)}{\lambda_{\star}\bar{F_{\star}}(x)+\lambda_{\dagger} \bar{F_{\dagger}}(x)}R_{\star}(dx)\\
&=\int_0^\infty\frac{\lambda_{\dagger} \psi(e^{-y})}{\lambda_{\star} e^{-y}+\lambda_{\dagger} \psi(e^{-y})}dy=\int_0^1 \frac{\lambda_{\dagger} \psi(u)}{\lambda_{\star} u+\lambda_{\dagger} \psi(u)} \frac{du}{u},
\end{split}
\end{equation*}
where the second equality is quite clear if $F_{\star}$ is strictly increasing, but holds in general (see e.g.~the \cite{FT12} for a discussion and justification).

If $\phi(\infty)<\infty$, then dominated convergence yields
\begin{equation*}
\lim_{\varepsilon\downarrow0}\int_\varepsilon^{2\varepsilon}\frac{\lambda_{\dagger}\psi(u)}{(\lambda_{\star}u+\lambda_{\dagger}\psi(u))u}du
=\lim_{\varepsilon\downarrow0}\int_0^1\frac{\lambda_{\dagger}\psi(u)}{(\lambda_{\star}u+\lambda_{\dagger}\psi(u))u}I_{(\varepsilon,2\varepsilon)}(u)du=0.
\end{equation*}

Since
\begin{equation*}
0\leq\int_\varepsilon^{2\varepsilon}\frac{\lambda_{\dagger}\psi(\varepsilon)}{(\lambda_{\star}u+\lambda_{\dagger}\psi(\varepsilon))u}du\leq\int_\varepsilon^{2\varepsilon}\frac{\lambda_{\dagger}\psi(u)}{(\lambda_{\star} u+\lambda_{\dagger} \psi(u))u}du,
\end{equation*}
we get that
\begin{equation}\label{lim1}
\lim_{\varepsilon\downarrow0}\int_\varepsilon^{2\varepsilon}\frac{\lambda_{\dagger}\psi(\varepsilon)}{(\lambda_{\star}u+\lambda_{\dagger}\psi(\varepsilon))u}du=0.
\end{equation}

Solving the integral,

\begin{equation*}
	\int_\varepsilon^{2\varepsilon}\frac{\lambda_{\dagger}\psi(\varepsilon)}{(\lambda_{\star}u+\lambda_{\dagger}\psi(\varepsilon))u}du
	=\log\bigg(1+\frac{\lambda_{\dagger}\psi(\varepsilon)}{\lambda_{\star}2\varepsilon+\lambda_{\dagger}\psi(\varepsilon)}\bigg),
\end{equation*}
and from  (\ref{lim1})
\begin{equation*}
	\lim_{\varepsilon\downarrow0}\log\bigg(1+\frac{\lambda_{\dagger}\psi(\varepsilon)}{\lambda_{\star}2\varepsilon+\lambda_{\dagger}\psi(\varepsilon)}\bigg)=0\Rightarrow \lim_{\varepsilon\downarrow0}\frac{\lambda_{\dagger}\psi(\varepsilon)}{\lambda_{\star}2\varepsilon+\lambda_{\dagger}\psi(\varepsilon)}=0\Rightarrow\lim_{\varepsilon\downarrow0}\frac{\lambda_{\dagger}\psi(\varepsilon)}{\lambda_{\star}\varepsilon}=0,
\end{equation*}

We may thus find
$\delta\in(0,1)$ such that $\lambda_{\dagger}\psi(u)\leq\lambda_{\star} u$ for each  $u\in (0,\delta)$. Thus,
\begin{equation*}
	\int_0^\delta\frac{\lambda_{\dagger}\psi(u)}{\lambda_{\star}u}\frac{du}{u}
	\leq
	2\int_0^\delta\frac{\lambda_{\dagger}\psi(u)}{\lambda_{\star} u+\lambda_{\dagger}\psi(u)}\frac{du}{u}\\
	\leq 2\int_0^1\frac{\lambda_{\dagger}\psi(u)}{\lambda_{\star} u+\lambda_{\dagger}\psi(u)}\frac{du}{u}<\infty.
\end{equation*}

Thus, since also $\psi(u)\leq1$, 
we have
\begin{equation*}
\frac{\lambda_{\dagger}}{\lambda_{\star}}\int_0^1\frac{\psi(u)}{u^2}du=\frac{\lambda_{\dagger}}{\lambda_{\star}}\int_0^\delta\frac{\psi(u)}{u^2}du+\frac{\lambda_{\dagger}}{\lambda_{\star}}\int_\delta^1\frac{\psi(u)}{u^2}du<\infty.
\end{equation*}

Rewriting (\ref{eq:expM}) in terms of $\psi$, we have
\begin{equation}\label{eq:em}
\mathbb{E}[M]=\frac{\lambda_{\dagger}}{\lambda_{\star}}\int_0^\infty\frac{\bar{F_{\dagger}}(x)}{\bar{F_{\star}}(x)}R_{\star}(dx)=\frac{\lambda_{\dagger}}{\lambda_{\star}}\int_0^\infty\frac{\psi(e^{-y})}{e^{-y}}dy=\frac{\lambda_{\dagger}}{\lambda_{\star}}\int_0^1\frac{\psi(u)}{u^2}du,
\end{equation}
and the finiteness of the latter expression establishes that (i)$\Rightarrow$(iii).

(iii)$\Leftrightarrow$(iv) follows readily from~(\ref{eq:em}).
\fim

\noindent\textbf{\textit{Proof of Theorem \ref{teo1}.}}
By Proposition \ref{prop:equiv}, if $\int_0^\infty\frac{\bar{F_{\dagger}}(x)}{\bar{F_{\star}}(x)}R_{\star}(dx)<\infty$, then
$M=\#(D\cap\hat{\Pi}_{\dagger})<\infty$ a.s.
Given $\{M<\infty\}$, consider the index of the last event of extinction  in $D$, namely $\hat{j}\doteq\max\{j:(S_j,Y_j)\in D\}$,
and note that there is only one index $k\geq1$ such that $(S_{\hat{j}},Y_{\hat{j}})\in D_k$, namely $\hat{k}$.
The next fitness record will happen at time  $T_{I_{\hat{k}+1}}$; then
	$\eta_0(s)\neq\varnothing$ for each $s\geq T_{I_{\hat{k}+1}}$,
since	 
$D\cap\hat{\Pi}_{\dagger}\cap[T_{I_{\hat{k}+1}},\infty)\times[0,\infty)=\varnothing$.
Hence,
$\{s\geq0:\eta_0(s)=\varnothing\}\subset[0,T_{I_{\hat{k}+1}})$
and 
$\{s\geq0:\eta_0(s)=\varnothing\}$ is a.s.~a bounded set.

If 
    $\int_0^\infty\frac{\bar{F_{\dagger}}(x)}{\bar{F_{\star}}(x)} R_{\star}(dx)=\infty$, then
	$M=\#(D\cap\hat{\Pi}_{\dagger})=\infty$ a.s.
	Given that $M=\infty$, suppose by contradiction that there is $\hat{s}>0$ such that
$\{s\geq0:\eta_0(s)=\varnothing\}\subset[0,\hat{s}].$
Since  $\mu_{\dagger}([0,\hat{s})\times[0,\infty))=\lambda_{\dagger}\hat{s}<\infty$, we have that
$\#([0,\hat{s})\times[0,\infty)\cap\hat{\Pi}_{\dagger})<\infty$ a.s.
Thus
$\#([0,\hat{s})\times[0,\infty)\cap D\cap\hat{\Pi}_{\dagger})<\infty$ {a.s.},
and 
	$\#([\hat{s},\infty)\times[0,\infty)\cap D\cap\hat{\Pi}_{\dagger})=\infty$ a.s. 
Now choose an extinction event $(S_j,Y_j)\in [\hat{s},\infty)\times[0,\infty)\cap D$ and observe that 
$\eta_0(S_j)=\varnothing$. The contradiction shows that $\{s\geq0:\eta_0(s)=\varnothing\}$ is a.s.~an unbounded set.
\fim

\noindent\textbf{\textit{Proof of Proposition \ref{marc}.}} The proof is analogous to the proof of Proposition \ref{marc:lambda} with the marking given by the kernel $K_{\dagger}$, replacing $(\lambda_{\star},{\bar{F}}_{\star})$ by $(\lambda_{\dagger},{\bar{F}}_{\dagger})$.
\fim

\noindent\textbf{\textit{Proof of Proposition \ref{marc:lamb}.}} The proof is analogous to  Proposition \ref{marc:mu}, with the independent $\mathbb{Q}_{\star}$-marking given by the measure induced by $X$. \fim

\noindent\textbf{\textit{Proof of Theorem \ref{teo2}.}}
By the homogeneity of the model, the process $\eta_t=\{ \eta_t(t+s): s\geq0\}$ has the same distribution for each $t \in \mathbb{R}$. 
In particular, 
\begin{equation}\label{dist}
\mathbb{P}(\eta_0(t)\in \cdot\, |\eta_0(0)=A)=\mathbb{P}(\eta_{-t}(0)\in \cdot\, |\eta_{-t}(-t)=A),
\end{equation} 
where $A$ is a locally finite subset of $[0,\infty)$. 
Let $(t_l)_{l\geq1}$ be an increasing sequence in ${\mathbb{R}}^+$ such that $t_l \rightarrow\infty$ as $l \rightarrow\infty$ and $A$ as above; take a sequence of processes $(\eta_{-t_l})_{l\geq1}$ such that 
$\eta_{-t_l}(-t_l)=A$, $l\geq1.$
Define the sequence $(\eta'_l)_{l\geq1}$ by 
\begin{equation*}
\eta'_l\doteq
\begin{cases}
\{X_i:T_i\in(-t_l,0]\}, \qquad\qquad\qquad\qquad\qquad\qquad\qquad \text{if} \quad \Pi_{\dagger}\cap(-t_l,0)=\varnothing\\
\{X_i:T_i \in (S_{J_1},0]\}\cup\Big(\bigcup_{k=1}^{\hat{k}_l-1}\{X_i>Y_{J_k}:T_i \in (S_{J_{k+1}},S_{J_k}]\}\Big)\\
\qquad\qquad \cup \{X_i>Y_{J_{\hat{k}_l}}:T_i \in (-t_l,S_{J_{\hat{k}_l}}]\}, \qquad  \text{otherwise},
\end{cases}
\end{equation*}
where 
$\hat{k}_l\doteq\max\{k: S_{J_k}>-t_l\},$
and  the union $\cup_{k=1}^{\hat{k}_l-1}$ above is empty if $\hat{k}_l=1$.
Note that $(\hat{k}_l)_{l\geq1}$ is an increasing sequence with $\hat{k}_l\rightarrow\infty$ as $l\rightarrow\infty$.  Thus, the sequence $(\eta'_l)_{l\geq1}$ is increasing and 
$\lim_{l\rightarrow\infty}\eta'_l=\bigcup_{l\geq1}\eta'_l=\hat{\eta}$,
with 
 \begin{equation*}\label{lim}
\hat{\eta}\doteq\{X_i:T_i \in (S_{J_1},0]\}\cup\Big(\bigcup_{k\geq1}\{X_i>Y_{J_k}:T_i \in (S_{J_{k+1}},S_{J_k}]\}\Big).
\end{equation*}

Define the sequence $(\eta''_l)_{l\geq1}$  by
\begin{equation*}
\eta''_l\doteq
\begin{cases}
A, &\text{if}\quad \Pi_{\dagger}\cap(-t_l,0)=\varnothing\\
A\cap(Y_{J_{\hat{k}_l}},\infty), &\text{otherwise}.
\end{cases}\\
\end{equation*}
From Proposition \ref{rec} for the random variables $(Y_j)$, 
${Y_{J_k}\rightarrow\infty}$ as $k\rightarrow\infty$, since $R_{\dagger}(x)\rightarrow\infty$ as $x\rightarrow\infty$. Thus, $(Y_{J_{\hat{k}_l}})_{l\geq1}$ is an increasing sequence with $\lim_{l\rightarrow\infty}Y_{J_{\hat{k}_l}}=\infty$, since $\hat{k}_l\rightarrow\infty$ as $l\rightarrow\infty$. Therefore, $(\eta''_l)_{l\geq1}$ is a decreasing sequence and
$\lim_{l\rightarrow\infty}\eta''_l=\bigcap_{l\geq1}\eta''_l=\varnothing.$
Because $\eta_{-t_l}(0)=\eta'_l\cup\eta''_l$ for each $l\geq1$,
we have that
$\lim_{l\rightarrow\infty}\eta_{-t_l}(0)=\lim_{l\rightarrow\infty}\eta'_l\cup\lim_{l\rightarrow\infty}\eta''_l=\hat{\eta},$
and notice that
the limit does not depend on the choice of $(t_l)_{l\geq1}$. Thus
$\lim_{t\rightarrow\infty}\eta_{-t}(0)=\hat{\eta}$ a.s.,
and, by (\ref{dist}), $\eta_0(t)\rightarrow\hat{\eta}$ in distribution as $t\rightarrow\infty.$
\fim

\noindent\textbf{\textit{Proof of Proposition \ref{equiv2}.}} The proof is analogous to that of Proposition \ref{prop:equiv}, using the same ideas, and changing the roles of $(\lambda_{\star},{\bar{F}}_{\star})$ and $(\lambda_{\dagger},{\bar{F}}_{\dagger})$.
\fim

\noindent\textbf{\textit{Proof of Theorem \ref{teo3}.}}
$N_0$ has a Poisson distribution with parameter $\Sigma_0$. 
Since 
$\Sigma_0=\mu_{\star}(E_0)=\lambda_{\star}(0-S_{-1})\sim$ Exponential($\lambda_{\dagger}/{\lambda_{\star}}$),
we have that
$\mathbb{E}[N_0]=
\mathbb{E}[\Sigma_0]=\lambda_{\star}/\lambda_{\dagger}<\infty,$
and thus $N_0<\infty$ a.s.

By (\ref{numlim}) and Proposition \ref{equiv2}, if 
$\int_0^\infty\frac{\bar{F_{\star}}(x)}{\bar{F_{\dagger}}(x)}R_{\dagger}(dx)<\infty,$ then $N<\infty$ a.s., and thus $\#\hat{\eta}<\infty$ a.s.
Otherwise, $N=\infty$ a.s., and we have $\#\hat{\eta}=\infty$ a.s.
$\hfill\square$

\end{document}